\documentclass[10pt,a4paper,twoside, reqno, tbtags]{amsart}

\usepackage{amsfonts, amssymb, amscd, amsmath, enumerate, verbatim}

\newcommand{\filebegin}{\begin{document}}
\newcommand{\fileend}{\end{document}}
\makeatother
\def\thefootnote{}
\newcommand{\lo}{\longrightarrow}
\newcommand{\NMM}{\hspace*{2mm}}

\renewcommand{\baselinestretch}{1.1}
\def\n{\noindent}

\numberwithin{equation}{section}
\def\mapdown#1{\Big\downarrow\rlap
{$\vcenter{\hbox{$\scriptstyle#1$}}$}}
\newtheorem{thm}{Theorem}[section]
\newtheorem{lem}{Lemma}[section]
\newtheorem{cor}{Corollary}[section]
\newtheorem{prop}{Proposition}[section]

\theoremstyle{definition}

\newtheorem{Defn}{Definition}[section]
\newtheorem{remark}{Remark}[section]
\newtheorem{result}{Result}[section]
\newtheorem{note}{Note}[section]
\newtheorem{notations}{Notations}
\newtheorem{Exmp}{Example}[section]

\begin{document}

\setcounter{page}{1} \noindent

\vspace*{2cm}
\begin{center}
{\bf Additivity of local function and dynamical system}
 \\[0.5cm]
{Sk. Selim$^{a}$, Chhapikul Miah$^{b}$, Monoj Kumar Das$^{c}$ and Shyamapada Modak$^{d*}$,  \footnote{$^*$Corresponding Author} \\[2mm]
$^{a}$Former Research Scholar, Department of Mathematics, University of Gour Banga, India\\[2mm]
$^{b}$Sukanta Mahavidyalaya, Dhupguri, Jalpaiguri 735210, West Bengal,
India\\[2mm]
$^{c,d}$Department of  Mathematics, University of Gour Banga\\Mokdumpur, Malda 732103, West Bengal, India \\[2mm]
{\tt E-mail: skselim2012@gmail.com}\\
{\tt E-mail: chhapikul.77@gmail.com}\\
{\tt E-mail: dmonojkr1@gmail.com}\\
{\tt E-mail: spmodak2000@yahoo.co.in}
} \\[2mm]
\end{center}%
\vspace*{0.5cm}
\begin{quotation}
\noindent
{\footnotesize
{\sc Abstract.}
The study of local function in topological spaces is  remarkable. Various branches have been developed through this study.  In this paper, we further consider the local function and exploring the various properties of the same by considering some generalized open sets. In this situation some of the properties of local function fails to hold due to the finite intersection property of the topology.  Due to this outcome, we are   investigating the situation of dynamical system and Topological Transitivity.
}
\end{quotation}
\\

\n \textbf{2020 Mathematics subject classification: }
Primary: 54H20; Secondary: 54G05, 54A10, 54C05, 54C10.\\

\markboth
{}
 {}



\section{\bf Introduction }\label{sec1}
In the field of pure and applied mathematics as well as in physics and chemistry, the study of limit points is an crucial part. In the study of mathematics this has been redefined through the Jankovic and Hamlett's \cite{JD1999} modern notation, by the name of local function. It is a part of an ideal on the topological spaces and the ideal was introduced by Kuratowski \cite{KK1966} and Vaidyanathswamy \cite{RV1960}. The ideal is a mathematical structure of subsets of a given set and satisfying the hereditary and the finite additive property. Many mathematicians \cite {A2021, TRD1990, JH1990, SC2007, SM2012, PNB2020} have given their reflection in the field of local function.

In this research article we are exploring the coincidence of local function of set and the limit point of a set. To do this, we will establish the relation between various local functions of literature with the limit points and to give counter examples. The remarkable discussion of this paper is, how one can  goes to limit point of a set via local function. That is numerically, we will approximating the limit points of a set through one by one local function.

\section{\bf Preliminaries}

A subset $A$ of a topological space $\mathbb{T}$ (where $T$ is a nonempty set and $\tau_T$ is a topology on $T$) is said to semi-open \cite{L1963} (resp. preopen \cite{MEE}, $b$-open \cite{DA}, $\beta$-open \cite{EEM} (or semi preopen \cite{AD}) if $A\subseteq Cl(Int(A))$ (resp. $A\subseteq Int(Cl(A))$, $A\subseteq Int(Cl(A)) \cup Cl(Int(A))$, $A\subseteq Cl (Int(Cl(A)))$), where `$Int$' and `$Cl$' stand for the `Interior' and `Closure' operators respectively. $\mathcal{O}(t)$ is the collection of all open sets containing $t$ where as $\mathcal{O}^s (t)$ stands for the collection of all semi-open sets containing $t$. Again $\mathcal{O}^p (t)$ (resp. $\mathcal{O}^{b}(t)$, $\mathcal{O}^{\beta} (t)$) denotes for the collection of all preopen (resp. $b$-open, $\beta$-open) sets containing $t$ and the collection of all $b$-open sets containing $t$. The `semi-closure', `pre-closure', `$b$-closure' and `$\beta$-closure' operators will be denoted as `$sCl$', `$pCl$' `$bCl$' and `$\beta Cl$' respectively. They are defined similar fashion with closure operator in the topological space.

We have already mentioned the idea of an ideal. Now, the local function of a set $A$ in an ideal $\mathbb{I}_T$ related  topological space \cite{JD1999} $(\mathbb{T}, \tau_T)$ ( or simply called $\mathbb{I}_T$ ideal on a topological space $\mathbb{T}$) is, $A^{\star} = \{t\in T| \; U\cap A \notin \mathbb{I}_T,\; U \in \mathcal{O}(t) \}$. This local function induces a finer topology on $T$ and its closure operator is defined as, $Cl^{\star}(A) = A\cup A^{\star}$. A subset $\bf A$ of a topological space $\mathbb{T}$ with an ideal $\mathbb{I}_T$ is called $\mathbb{I}$-dense if and only if ${\bf A}^{*}=\mathbb{T}$ \cite{DGR1999}. Equivalently ${\bf A}$ is $\mathbb{I}_T$-dense in $\mathbb{T}$ if and only if for every non-empty open subset ${\bf O}$ of $\mathbb{T}$, ${\bf A}\cap {\bf O}\notin \mathbb{I}_T$. Every $\mathbb{I}_T$-dense subset of a topological space is a dense set \cite{RV1960}.

Before consideration the other local functions and its related consequences, we are trying to recognising the local function in new format.

\begin{Exmp}
\label{Example 1}
Let $T = \mathbb{R}$, set of reals, $\tau_{T} = \{\emptyset, \mathbb{R}, \mathbb{Q} \},\; \mathbb{I}_T = \wp(\mathbb{Q})$, where  $\mathbb{Q}$ and $\wp(\mathbb{Q})$ denote the set of rationals and the set of all subsets of $\mathbb{Q}$ respectively. Then one of a basis of the $\star$ topology (or $\tau_{T}^{\star}$ topology) is, $\mathbb{B}(\mathbb{I}, \tau_{T}) = \{\mathbb{R}\} \cup \wp(\mathbb{Q}) \cup \{ (\mathbb{R} \setminus \mathbb{Q}) \cup A|\; A\in \wp(\mathbb{Q}) \}$. Then $1 \in \mathbb{Q}^{\prime}(\tau_{T})$, derived set w.r.to $\tau_{T}$ topology but $1 \notin \mathbb{Q}^{\prime}(\tau_{T}^{\star})$, derived set w.r.to $\tau_{T}^{\star}$ topology. Whereas $\mathbb{Q}^{\star} = \emptyset$. Furthermore, for $i \in \mathbb{R} \setminus \mathbb{Q}$,  $(\{i \})^{\star} = \mathbb{R} \setminus \mathbb{Q}$. But, $(\{i \})^{\prime}(\tau_{T}) = (\mathbb{R}\setminus \mathbb{Q}) \setminus \{i\} = (\{i \})^{\prime}(\tau_{T}^{\star})$.
\end{Exmp}

Due to this example, the reverse inclusion of $A^{\star} \subseteq A^{\prime}(\tau_{T})$, $A^{\star} \subseteq A^{\prime}(\tau_{T}^{\star})$ (mentioned that local function of a set w.r to the topology $\tau_{T}$ and the $\star$ topology are same) and $A^{\prime}(\tau_{T}^{\star}) \subseteq A^{\prime}(\tau_{T})$ are not true. But interestingly the operator $Cl^{\star}: \wp(T) \mapsto \wp(T)$ (defined above) induces a topology on $T$.

 Local function related to generalized open sets are:

$A^{\star s}$ \cite{KN2010, SS2012} (resp. $A^{\star p}$ \cite{MM2012, A2021}, $A^{\star b}$ \cite{TB2017}, $A^{\star \beta}$ \cite{M2011, TB2017}) $ = \{t \in T|\; U\cap A \notin \mathbb{I}_T \}$, where $U \in \mathcal{O}^s(t)$ (resp. $U \in \mathcal{O}^{p}(t)$, $U \in \mathcal{O}^b(t)$, $U \in \mathcal{O}^{\beta}(t)$).

Closure operator related local function is,  $\Gamma (A)$ \cite{AN2013, KOAM} (resp. $\gamma (A)$ \cite{IM2018}, $\xi_s(A)$ \cite{YK2019}, $\xi_p(A)$, $\xi_b(A)$, $\xi_{\beta}(A)$ \cite{BTP2023}) $ = \{z\in T |\; Cl(U^o) \cap A \notin \mathbb{I}_T \}$ (resp. $\{z\in T |;\ sCl(U^o) \cap A \notin \mathbb{I}_T \}$, $\{z\in T ;\ sCl(U^s) \cap A \notin \mathbb{I}_T \}$, $\{z\in T ;\ pCl(U^p) \cap A \notin \mathbb{I}_T \}$, $\{z\in T ;\ bCl(U^b) \cap A \notin \mathbb{I}_T \}$, $\{z\in T |;\ \beta Cl(U^{\beta}) \cap A \notin \mathbb{I}_T \}$), where $U^o \in \mathcal{O}(z)$, $U^s \in \mathcal{O}^{s}(z)$, $U^p \in \mathcal{O}^{p}(z)$, $U^b \in \mathcal{O}^{b}(z)$, $U^{\beta} \in \mathcal{O}^{\beta}(z)$).

 The associated set-valued set function\cite{KOAM, MS2021} of the above operators are as follows: $\psi_{\xi_s}(A) $\cite{TN1986, YK2021}$ = T\setminus \xi_s(T\setminus A)$;  $\psi_{\xi_p}(A) = T\setminus \xi_p(T\setminus A)$;  $\psi_{\xi_b}(A) = T\setminus \xi_b(T\setminus A)$;  $\psi_{\xi_{\beta}}(A) $\cite{BTP2023}$ = T\setminus \xi_{\beta}(T\setminus A)$; $\psi_{\Gamma}(A) = T\setminus \Gamma(T\setminus A)$ \cite{AN2013}.

\section{\bf  Coincident of local functions}

 In this section, we shall try to give more answers of Yalaz and Kaymakcı's \cite{YK2023} question: ``{\bf Question 1}: For a subset $M$ in any $I$-space, are local function and $\xi_{\Gamma}^*$-local function always comparable with respect to the subset relation? So is it always either $M^* \subseteq \xi_{\Gamma}^*(M)$ or $\xi_{\Gamma}^*(M) \subseteq M^*$?"

For a subset $A$ of an ideal topological space $(T,\tau_{T},\mathbb{I}_T)$, $\xi^{\star}$ local function \cite{YK2023} is defined as, $\xi^{\star}(A) = \{t\in T|\; O^{\star} \cap A \notin \mathbb{I}_{T} \}$, where $O \in \mathcal{O}(t)$. Its associated set-valued set function is $\psi_{\xi^{\star}}(A) = T\setminus \xi^{\star}(T\setminus A)$ \cite{YK2023}.

In Theorem 3.5 of \cite{YK2023}, the authors have discussed equality of three local functions. But for the equality of these, they have used the condition $\mathbb{I}_{T} \cap \tau_{T} = \{\emptyset \}$. Because for $\mathbb{I}_{nw}$, $\mathbb{I}_{nw}$ does not contain any nonempty open set (where $\mathbb{I}_{nw}$ denotes the collection of all nowhere dense sets in the topological space $\mathbb{T}$). If $\emptyset \neq O \in \mathbb{I}_{nw}$, then $Int (Cl(O)) = \emptyset$, which is a contradiction. Followings are the more general answer of the above question.

\begin{lem} \cite{SC2006}
\label{ }
Let $\mathbb{I}_{T}$ be an ideal on a topological space $\mathbb{T}$. Then $\mathbb{I}_{T} \cap \tau_{T} = \{\emptyset \}$ if and only if, for any nonempty open set $O$, $O^{\star} = Cl(O)$.
\end{lem}

\begin{thm}
\label{ }
Let $\mathbb{I}_{T}$ be an ideal on a topological space $\mathbb{T}$ such that $\mathbb{I}_{T} \cap \tau_{T} = \{\emptyset \}$. Then $\xi^{\star}(A) = \Gamma (A)$ for all $A\in \wp(T)$.
\end{thm}

Conversely, we can get the following,

\begin{cor}
\label{ }
Let $\mathbb{I}_{T}$ be an ideal on a topological space $\mathbb{T}$.  If $\xi^{\star}(A) = \Gamma (A)$ for all $A\in \wp(T)$, then $\mathbb{I}_{T} \cap \tau_{T} = \{\emptyset \}$.
\end{cor}

\begin{proof}
Suppose $O$ be a nonempty open set in $T$. Given that
 $\xi^{\star}(O) = \Gamma (O)$. Thus $O^{\star} \cap O \notin \mathbb{I}_{T}$ as well as $Cl(O) \cap O \notin \mathbb{I}_{T}$. Therefore $O \notin \mathbb{I}_{T}$. As $O$ is arbitrary, then $\mathbb{I}_{T} \cap \tau = \{\emptyset \}$.
\end{proof}

In these respect the Theorem 3.6 of \cite{YK2023} become,

\begin{thm}
\label{ }
Let $\mathbb{I}_{T}$ be an ideal on a topological space $\mathbb{T}$. Then for any $A\in \wp(T)$, $A^{\star} = \xi^{\star}(A) = \Gamma(A)$ if and only if $\mathbb{I}_{T}\cap \tau_{T} = \{\emptyset \}$.
\end{thm}

As the consequence of the above results of this section we get,

\begin{lem}
\label{ }
Let $\mathbb{I}_{T}$ be an ideal on a topological space $\mathbb{T}$. Then for any $A\in \wp(T)$, $\psi_{\xi^{\star}}(A) = \psi_{\Gamma}(A)$, if and only if $\mathbb{I}_{T} \cap \tau_{T} = \{\emptyset \}$.

\end{lem}

Due to this lemma we get the following coincident results:

For an ideal $\mathbb{I}_{T}$ on a topological space $\mathbb{T}$ with $\mathbb{I}_{T} \cap \tau_{T} = \{\emptyset \}$.

\begin{thm}
\label{ }
\begin{enumerate}
  \item The topology $\{ M\subseteq T| \; M \subseteq \psi_{\xi^{\star}}(M) \}$ \cite{YK2023} on $T$ and the topology $\{ M\subseteq T| \; M \subseteq \psi_{\Gamma}(M) \}$ \cite{AN2013} on $T$ are equal.
  \item The topology $\{ M\subseteq T| \; M \subseteq Int (Cl( \psi_{\xi^{\star}}(M))) \}$ \cite{YK2023} on $T$ and the topology $\{ M\subseteq T| \; M \subseteq Int(Cl( \psi_{\Gamma}(M))) \}$ \cite{AN2013} on $T$ are equal.
\end{enumerate}
\end{thm}

\section{\bf Local function apart from coincident}
In this section, we consider new type of local functions using $O \in \mathcal{O}^k(t)$ sets where and $k \in \{s,\;p,\; b,\; \beta\}$.

Now, we consider an example in support of ``$\mathbb{I}_{T} \cap \tau_{T} = \{\emptyset \}$" does not mean that $A^{\star k} = Cl^k(A)$.

\begin{Exmp}
\label{Example 2}
$(i)$ Let $T = \{t^1, t^2, t^3, t^4 \},\; \tau_T = \{\emptyset, T, \{t^3, t^4 \}, \{t^2, t^3, t^4 \}, \{t^1, t^3, t^4 \} \},\; \\ \mathbb{I}_{T} = \{\emptyset, t^1 \}$. Let $A = \{t^1, t^3 \}$, then $A^{\star p} = \{t^3 \}$. But $Cl^p(A) = \{t^2, t^4 \}$. Further, $A^{\star \beta} = \{t^3 \}$, but $Cl^{\beta} (A) = \{t^2, t^4 \}$, whereas $\tau_T \cap \mathbb{I}_{T} = \{ \emptyset \}$.

$(ii)$ Let $T = \{t^1, t^2, t^3, t^4 \},\; \tau_{T} = \{\emptyset, T, \{t^1 \}, \{t^2 \}, \{t^1, t^2 \} \},\;  \mathbb{I}_{T} = \{\emptyset, t^3 \}$. Then, for $A = \{t^1, t^3 \}$, $A^{\star s} = \{t^3 \}$. But $Cl^s(A) = \{t^2, t^4 \}$, whereas $\tau_{T} \cap \mathbb{I}_{T} = \{\emptyset \}$.
\end{Exmp}

Following local functions are the new types of local function:

\begin{Defn}
\label{Definition 1}
Let $\mathbb{I}_{T}$ be an ideal on a topological space $\mathbb{T}$ and $A\subseteq T$. Then $k$  type local function is defined as, $A^{\diamond k} = \{t\in T |\; O^{\star k} \cap A \notin \mathbb{I} \}$, where $O \in \mathcal{O}^k(t)$ and $k \in \{s,\;p,\; b,\; \beta\}$.
\end{Defn}

\begin{thm}
\label{ }
Let $\mathbb{I}_T$ and $\mathbb{J}_{T}$ be two ideals on a topological space $\mathbb{T}$ and $A, B \subseteq T$. Then for $k\in \{s,\;p,\; \beta\}$,
\begin{enumerate}
  \item $A^{\diamond k} \subseteq B^{\diamond k}$, when $A\subseteq B$.
  \item $A^{\diamond k}(\mathbb{J}_{T}) \subseteq A^{\diamond k}(\mathbb{I}_{T})$, when $\mathbb{I}_{T} \subseteq \mathbb{J}_{T}$.
  \item $A^{\diamond k} \subseteq Cl^k(A)$.
  \item $I^{\diamond k} = \emptyset$, when $I \in \mathbb{I}_{T}$.
  \item $\emptyset^{\diamond k} = \emptyset$.
\end{enumerate}
\end{thm}

Interrelations between various local functions are follows after the following examples:

\begin{Exmp}
\label{Example 3}
\begin{enumerate}
  \item Let $T = \{t_1, t_2, t_3 \},\; \tau_{T} = \{\emptyset, T, \{t_1 \} \},\; \mathbb{I}_{T} = \{\emptyset, \{t_3 \} \}$. Then, $PO(\mathbb{T})$ (collection of all  preopen sets)$=\{\emptyset, T, \{t_1 \}, $ $ \{t_1, t_3 \}, \{t_1, t_2 \} \}, \emptyset^{\star p} = \emptyset, \; T^{\star p } = T,\; (\{t_1, t_3 \})^{\star p} = T= (\{t_1, t_2 \})^{\star p}$ and $(\{t_2, t_3 \})^{\star p} = \{t_2 \}, (\{t_2, t_3 \})^{\diamond p} =T$.
  \item  Let $T = \{t_1, t_2, t_3, t_4 \},\; \tau_{T} = \{\emptyset, T, \{t_3, t_4 \}, \{t_2, t_3, t_4 \}, \{t_1, t_3, t_4 \} \},\; \mathbb{I}_{T} = \{\emptyset, \{t_1 \} \}$. Then $PO(\mathbb{T}) = \{\emptyset, T, \{t_3 \}, \{t_4 \}, \{t_3, t_4 \}, \{t_1, t_3 \}, \{t_2, t_3 \}, \\ \{t_1, t_4 \}, \{t_2, t_4 \}, \{t_1, t_2, t_3 \}, \{t_1, t_2, t_4 \}, \{t_2, t_3, t_4 \}, \{t_1, t_3, t_4 \} \}, \\ \emptyset^{\star p} = \emptyset,\; T^{\star p} = T = (\{t_3, t_4 \})^{\star p} = (\{t_2, t_3, t_4 \})^{\star p} = (\{t_1, t_3, t_4 \})^{\star p},\; (\{t_3 \})^{\star p} = \{t_3 \},\; (\{t_4 \})^{\star p} = \{ t_4 \},\; (\{t_1, t_3 \})^{\star p} = \{t_3 \},\; (\{t_2, t_3 \})^{\star p} = \{t_2, t_3 \},\; (\{t_1, t_4 \})^{\star p} = \{t_4 \},\; (\{t_2, t_4 \})^{\star p} = \{t_2, t_4 \},  (\{t_1, t_2, t_3 \})^{\star p} =\{t_2, t_3 \},\; (\{t_1, t_2, t_4 \})^{\star p} = \{t_2, t_4 \}$.

     $(i)$ $ (\{t_1, t_3 \})^{\diamond p} = \{t_3 \},\; (\{t_1, t_4 \})^{\diamond p} = \{t_4 \},\; (\{t_1, t_3, t_4 \})^{\diamond p} = T$.

     $(ii)$ $pCl(\{t_1, t_2, t_3 \}) = \{t_1, t_2, t_3 \}$ whereas $(\{t_1, t_2, t_3 \})^{\star p} = \{t_2, t_3 \}$.

  \item Let $T = \{t_1, t_2, t_3 \},\; \tau_{T} = \{\emptyset, T, \{t_1, t_3 \} \},\; \mathbb{I}_{T} = \{\emptyset, \{t_1 \}, \{t_3 \}, \{t_1, t_3 \} \}$.
      Then,

      $SO(\mathbb{T}) (\text{collection of all semi-open sets}) =\{\emptyset, T, \{t_1, t_3 \} \},\; \emptyset^{\star s} = \emptyset, \; T^{\star s } = T = (\{t_1, t_3 \})^{\star s}, (\{ t_2, t_3 \})^{\star s} = \{t_2 \},\; (\{t_2, t_3 \})^{\diamond s} = T$.

  \item Let $T = \{t_1, t_2, t_3, t_4 \},\; \tau_{T} = \{\emptyset, T, \{t_1 \}, \{t_2 \}, \{t_1, t_2 \} \},\; \mathbb{I}_{T} = \{\emptyset, \{t_3 \} \}$. Then, $SO(\mathbb{T}) = \{\emptyset, T, \{t_1 \}, $ $\{t_2 \}, \{t_1, t_2 \}, \{t_1, t_3 \}, \{t_1. t_4 \}, \{t_2, t_3 \}, \{t_2, t_4 \}, \\ \{t_1, t_2, t_3 \}, \{t_1, t_2, t_4 \}, \{t_2, t_3, t_4 \}, \{t_1, t_3, t_4 \}\}$. Then, \\ $(\{t_1 \})^{\star s} = \{t_1 \},\; (\{t_1, t_2 \})^{\star s} = T = (\{t_1, t_2, t_3 \})^{\star s} = (\{t_1, t_2, t_4 \})^{\star s} = T^{\star s} = (\{t_1, t_2 \})^{\star s}, (\{t_1, t_3 \})^{\star s} = \{t_1 \},\; (\{t_1, t_4 \})^{\star s} = \{t_1, t_4 \},  (\{t_2 \})^{\star s} = \{t_2 \}, (\{t_2, t_3 \})^{\star s} = \{t_2 \},\; (\{t_2, t_4 \})^{\star s} = \{t_2, t_4 \} = (\{t_2, t_3, t_4 \})^{\star s},\;
       $ $(\{t_1, t_3, t_4 \})^{\star s} = \{t_1, t_4 \}$.

      $(i)$ $(\{t_1, t_3 \})^{\diamond s} = \{t_1 \},\; (\{t_2, t_3 \})^{\diamond s} = \{t_2 \},\; (\{t_1, t_2, t_3 \})^{\diamond s} = T$.

      $(ii)$ $sCl(\{t_2, t_3 \}) = \{t_2, t_3 \}$ whereas $(\{t_2, t_3 \})^{\star s} = \{t_2 \}$.

\end{enumerate}

\end{Exmp}

For a topological space $\mathbb{T}$, $\tau_{T} \subseteq SO(\mathbb{T}) (\mbox{resp.}\; PO(\mathbb{T}) ) \subseteq BO (\mathbb{T}) (\text{collection of all $b$-open sets}) \subseteq \beta O (\mathbb{T}) (\text{collection of all $\beta$-open sets}) $ holds. This relation establish the following theorem.

\begin{thm}
\label{ }
Let $\mathbb{I}$ be an ideal on a topological space $\mathbb{T}$. Then for $A\subseteq T$,
\begin{enumerate}
  \item $A^{\star \beta} \subseteq A^{\star b} \subseteq A^{\star s}\subseteq A^{\star}$.
  \item $A^{\star \beta} \subseteq A^{\star b} \subseteq A^{\star p}\subseteq A^{\star}$.
  \item $A^{\diamond \beta} \subseteq A^{\diamond b} \subseteq A^{\diamond s}\subseteq A^{\diamond}$.
  \item $A^{\diamond \beta} \subseteq A^{\diamond b} \subseteq A^{\diamond p}\subseteq A^{\diamond}$.
\end{enumerate}
\end{thm}

For proof of this theorem the reader may take the help of \cite{SIM2020}.

Example \ref{Example 3}(1,3) shows that $\diamond k$ local function and $\star k$ local function are independent to each other in respect of set inclusion relation. The Example \ref{Example 3}(2,4) also mentioned that $\diamond k$ local function is not additive.

\begin{thm}
\label{ }
Let $\mathbb{I}_{T}$ be an ideal on a topological space $\mathbb{T}$. Then for $A\subseteq T$, $(i)$ $A^{\diamond k} \nsubseteq kCl(A)$; $(ii)$ $A^{\star k} \subseteq \xi_{k}(A)$, where $k\in \{s,\;p,\; b$,\;$\beta\}$.
\end{thm}
$(i)$ follows from Example \ref{Example 3} (1) and $(ii)$ is obvious fact.

Equality of $A^{\diamond k}$ with $kCl(A)$ does not hold even if $\mathbb{I}_{T} \cap \tau_{T} = \{\emptyset \} = \mathbb{I}_{T} \cap KO (\mathbb{T})$, where $K = S, P, B$ and $\beta$. These have been followed by the Example \ref{Example 3} (2,4). However following hold:

\begin{cor}
\label{ }
Let $\mathbb{I}_{T}$ be an ideal on a topological space $\mathbb{T}$. Then for $O \in \mathcal{O}^{k}$, $O^{\star k} = kCl(O)$ implies $\mathbb{I} \cap \mathcal{O}^{k} = \{\emptyset \}$, where $k\in \{s,\;p,\; b$,\;$\beta\}$.
\end{cor}

Thus the $\diamond$ local functions ($()^{\diamond k}$) do not coincides with the following local functions:

$\xi_{k}^{\diamond }(A) = \{t\in T|\; O^{\star k} \cap A \notin \mathbb{I}_{T} \}$, where $O \in \mathcal{O}^{k}(t)$ and $k\in \{s,\;p,\; b$,\;$\beta\}$.

However $\xi^{\star}$ local function and $\Gamma$ local function are coincide each other when $\mathbb{I} \cap \tau_{T} = \{\emptyset \}$ (see Section 3).

Semi-local function and pre local function, these are independent each other with respect to the set inclusion relation. This has been followed by the following example:

\begin{Exmp}
\label{Example 4}
$(i)$ Let ${T} =\{t_1, t_2, t_3 \},\; \tau_{\mathbb{T}} = \{\emptyset, {T}, \{t_1, t_2 \} \},\; \mathbb{I}_{\mathbb{T}} = \{\emptyset, \{t_1 \}\}$. Then $SO(\mathbb{T}) = \{\emptyset, {T}, \{t_1, t_2 \} \},\\  PO(\mathbb{T}) = \{\emptyset, {T}, \{t_1, t_2 \}, \{t_1 \}, \{ t_2 \}, \{t_1, t_3 \}, \{t_2, t_3 \} \}$. Consider $A = \{t_2, t_3 \}$, then $A^{\star p} = \{t_2, t_3 \}$ and  $A^{\star s} = {T}$.

$(ii)$ Let ${T} =\{t_1, t_2, t_3, t_4 \},\; \tau_{\mathbb T} = \{\emptyset, {T},  \{t_1 \}, \{ t_3 \}, \{t_1, t_3 \}, \{t_1, t_3, t_4 \} \},\; \mathbb{I}_{\mathbb T} = \{\emptyset, \{t_2 \} \}$. Then  $SO(\mathbb{T}) = \{\emptyset, {T}, \{t_1 \}, \{ t_3 \}, \{t_1, t_3 \},  \{t_1, t_2 \}, \{t_1, t_4 \}, \{t_2, t_3 \},\\ \{t_3, t_4 \}, \{t_1, t_2, t_3 \}, \{t_1, t_2, t_4 \}, \{t_1, t_3, t_4 \}, \{t_2, t_3, t_4 \} \}, \\ PO(\mathbb{T}) = \{\emptyset, {T}, \{t_1 \}, \{ t_3 \}, \{t_1, t_3 \}, \{t_1, t_3, t_4 \}, \{t_1, t_2, t_3 \} \}$. Suppose $B = \{t_2, t_3, t_4 \}$, then $B^{*s} = \{t_3, t_4 \},  B^{*p} = \{t_2, t_3, t_4 \}$.
\end{Exmp}

For the decomposition of local function, we consider $\alpha$ set. In a topological space $\mathbb{T}$, a subset $A$ of $\mathbb{T}$ is mentioned as a $\alpha$ set \cite{NO1965} if it is semi-open as well as preopen set. That is $\tau_{\mathbb{T}}^{\alpha}$ (collection of all $\alpha$ sets)$ = SO(\mathbb{T}) \cap PO(\mathbb{T})$. Remarkable fact that $\tau_{T}^{\alpha}$ form a topology and the topology is finer than the original topology whereas $SO(\mathbb{T})$ and  $PO(\mathbb{T})$ do not form a topology in general.

\begin{thm}
\label{ }
Let $\mathbb{I}$ be an ideal on a topological space $\mathbb{T}$. Then for $A\subseteq T$, $(i)$ $A^{\star \beta} \subseteq A^{\star b} \subseteq A^{\star p}$ (resp. $A^{\star s}$) $\subseteq A^{\star \alpha} \subseteq A^{\star }$; $(ii)$ $A^{\diamond \beta} \subseteq A^{\diamond b} \subseteq A^{\diamond p}$ (resp. $A^{\diamond s}$) $\subseteq A^{\diamond \alpha} \subseteq A^{\diamond }$.
\end{thm}

Following example shows that the reverse inclusion of $A^{\star p}$ (resp. $A^{\star s}$) $\subseteq A^{\star \alpha}$ and $ A^{\diamond p}$ (resp. $A^{\diamond s}$) $\subseteq A^{\diamond \alpha}$ do not hold:

\begin{Exmp}
\label{Example 5}
$(a)$ Recall the Example \ref{Example 1}. Then $\mathcal{O}^s = \{\emptyset \} \cup \{A\supseteq \mathbb{Q}|\; A\subseteq \mathbb{R} \}$,  $\mathcal{O}^p = \{\mathbb{R} \cup \{A\subseteq \mathbb{R}|\;  A\nsubseteq \mathbb{R} \setminus \mathbb{Q} \}$ and $\tau_{T}^{\alpha} = \mathcal{O}^p \cap \mathcal{O}^s$. Suppose $i \in \mathbb{R} \setminus \mathbb{Q}$. Then $(\{i \})^{\star s} = \{i \}$ and $(\{i \})^{\star \alpha} = \mathbb{R} \setminus \mathbb{Q}$.

$(b)$ If we replace the ideal $\mathbb{I}_{T} = \wp(\mathbb{Q})$ with $\mathbb{I}_{T} = \{\emptyset, \{i \} \}$ (where $i$ is mentioned above) in the above example. Then for $p \in \mathbb{Q}$, $(\{i, p \})^{\star \alpha} = \mathbb{R}$ and  $(\{i, p \})^{\star p} = \mathbb{R} \setminus \mathbb{Q}$.
\end{Exmp}

\section{\bf Set with respect to ${\star sp}$ notation and transitivity}

At first we shall investigate the answer of the question, Is $A^{\star s} \cap A^{\star p} $ always nonempty?

\begin{thm}
\label{ }
Let $\mathbb{I}_T$ be an ideal on a topological space $\mathbb{T}$ and $A\subseteq T$. Then
\begin{enumerate}
\item for $\mathcal{O}^s \subseteq \mathcal{O}^p$ and $A^{\star s} \neq \emptyset \neq A^{\star p}$, $A^{\star s} \cap A^{\star p} \neq \emptyset$.
\item for $\mathcal{O}^p \subseteq \mathcal{O}^s$ and $A^{\star s} \neq \emptyset \neq A^{\star p}$, $A^{\star s} \cap A^{\star p} \neq \emptyset$.
\end{enumerate}
\end{thm}

\begin{proof}
1. Let $t\in A^{\star p}$. Then for all $U \in \mathcal{O}^p(t)$, $U\cap A \notin \mathbb{I}_T$. Since $\mathcal{O}^s \subseteq \mathcal{O}^p$, then for all $V \in \mathcal{O}^s(t)$, $V\cap A \notin \mathbb{I}_T$. Hence the result.

2. Similar with 1.
\end{proof}

This theorem gives a sufficient condition to hold $A^{\star s} \cap A^{\star p} \neq \emptyset$. But converse part of this theorem is an open question.

\begin{Defn}
\label{ }
Let $\mathbb{I}_T$ be an ideal on a topological space $\mathbb{T}$ and $A\subseteq T$. A point $t\in \mathbb{T}$ is called $\star sp$-derived point of $A$ if $t \in A^{\star s} \cap A^{\star p}$.

Collection of all $\star sp$-derived points of $A$ is denoted as $A^{\star sp}(\mathbb{I}_T)$ (or simply $A^{\star sp}$).
\end{Defn}
For the existence of $\star sp$-derived point, we consider the following example:
\begin{Exmp}
\label{Example 5.1.1}
Suppose $T=\{t_1, t_2, t_3\}$ is a set endowed with the topology $\tau=\{\emptyset, {T}, \{t_1\}, \{t_2\}, \{t_1, t_2\}\}$ and $\mathbb{I}_{T}=\{\emptyset,\{t_3\}\}$. Now $PO(\mathbb{T}) = \{\emptyset, {T}, \{t_1\}, \{t_2\}, \{t_1, t_2\}\}$ and $SO(\mathbb{T}) = \{\emptyset, {T}, \{t_1\}, \{t_2\}, \{t_1, t_2\}, $ $ \{t_2, t_3\}\}$. Now, $\{t_1\}^{*s}=\{t_1\}$, $\{t_1, t_2\}^{*s}={T}$, $\{t_2, t_3\}^{*s}=\{t_2, t_3\}$, $\{t_1, t_3\}^{*s}=\{t_1, t_3\}$, $\{t_2\}^{*s}=\{t_2, t_3\}$, $\{t_3\}^{*s}=\{t_3\}$, $\{\emptyset\}^{*s}=\{\emptyset\}$, $\{{T}\}^{*s}={T}$ and $\{t_1\}^{*p}=\{t_1, t_3\}$, $\{t_2\}^{*p}=\{t_2, t_3\}$, $\{t_3\}^{*p}=\emptyset$, $\{t_1, t_2\}^{*p}={T}$, $\{t_2, t_3\}^{*p}=\{t_2, t_3\}$, $\{t_3, t_1\}^{*p}=\{t_3, t_1\}$, $\{\emptyset\}^{*p}=\{\emptyset\}$, $\{{T}\}^{*p}={T}$.  Take, $A=\{t_1, t_2\}$, then $A^{\star sp}={T}$.
\end{Exmp}

\begin{thm}
\label{ }
Let $\mathbb{I}_T$ be an ideal on a topological space $\mathbb{T}$ and $A, B \in \wp(T)$. Then,
\begin{enumerate}
  \item $A^{\star sp} \subseteq B^{\star sp}$ if $A\subseteq B$.
  \item for an ideal $\mathbb{J}_T$ on $T$ with $\mathbb{J}_T \subseteq \mathbb{I}_T$, $A^{\star sp}(\mathbb{I}_T) \subseteq A^{\star sp}(\mathbb{J}_T)$.
  \item $A^{\star sp} \subseteq pCl(A)$, $A^{\star sp} \subseteq sCl(A)$.
  \item $(A^{\star sp})^{\star sp} \subseteq A^{\star sp}$.
  \item $A^{\star sp} \subseteq A^{\star s} \subseteq A^{\star \alpha} \subseteq A^{\star}$.
  \item $A^{\star sp} \subseteq A^{\star p} \subseteq A^{\star \alpha} \subseteq A^{\star}$.
  \item for ${\bf I} \in \mathbb{I}_T$, $(A\setminus {\bf I})^{\star sp} = (A\cup {\bf I})^{\star sp} = A^{\star sp}$.
  \item $A^{\star sp} \cup B^{\star sp}  \subseteq (A\cup B)^{\star sp}$.

\end{enumerate}
\end{thm}

Following example shows that $A^{\star sp} \cup B^{\star sp}  \neq (A\cup B)^{\star sp}$.

\begin{Exmp}
\label{Example 6}
Consider $A = \{t_1, t_4 \}$ in the Example \ref{Example 4}(ii). Then $A^{\star s} = \{t_1, t_4 \}$,  $A^{\star p} = \{t_1, t_2, t_4 \}$ and hence  $A^{\star sp} = \{t_1, t_4 \}$. Again, $B^{\star sp} = \{t_3, t_4 \}$ and hence $A^{\star sp} \cup B^{\star sp}= \{t_1, t_3, t_4 \}$. Thus, $(A\cup B)^{\star sp}=\{t_1, t_2, t_3, t_4 \} \neq A^{\star sp} \cup B^{\star sp}$.
\end{Exmp}

Due to this example, we conclude that $Cl^{\star sp}: \wp(T) \longmapsto \wp(T)$ defined by $Cl^{\star sp}(A) = A\cup A^{\star sp}$ is not a closure operator and hence not induces a topology on $T$.

\begin{lem}
\label{5.1}
Let $\mathbb{I}_{T}$ be an ideal on a topological space $\mathbb{T}$. Then, $\mathbb{I}_{T} \cap \mathcal{O}^s = \{\emptyset\}$ and $\mathbb{I}_{T} \cap \mathcal{O}^p = \{\emptyset\}$ if  for any nonempty set $A$, $A^{\star sp}$ is nonempty.
\end{lem}
\begin{proof}
Let $\emptyset \neq U\in \mathcal{O}^s$  and $\emptyset \neq V\in \mathcal{O}^p$. Since,  $U^{\star sp}$  and  $V^{\star sp}$ are nonempty, then $U=U\cap \mathbb{T}\notin  \mathbb{I}_{T}$ and $V=V\cap \mathbb{T}\notin  \mathbb{I}_{T}$. Since, $U\in \mathcal{O}^s$  and $V\in \mathcal{O}^p$ are arbitrary, then $\mathbb{I}_{T} \cap \mathcal{O}^s = \{\emptyset \}$ and $\mathbb{I}_{T} \cap \mathcal{O}^p = \{\emptyset \}$.
\end{proof}
The converse part may not be true. For this, we consider the following example:
 \begin{Exmp}
\label{Example 5.1.2}
In Example \ref{Example 5.1.1}, $\mathbb{I}_{T} \cap \mathcal{O}^s = \{\emptyset\}$ and $\mathbb{I}_{T} \cap \mathcal{O}^p = \{\emptyset\}$. If we consider $A=\{t_3\}$, then $A^{\star sp}=\emptyset$.
\end{Exmp}

\begin{prop}
\label{R5.1}
Let $\mathbb{I}_{T}$ be an ideal on a topological space $\mathbb{T}$. If $\mathbb{I}_{T} \cap \mathcal{O}^s = \{\emptyset \}$ and $\mathbb{I}_{T} \cap \mathcal{O}^p = \{\emptyset \}$, then $O^{\star sp}=sCl(O)\cap pCl(O)$ for all $O\in \mathcal{O}^s\cap \mathcal{O}^p$.
\end{prop}
\begin{proof}
Let $t\in O^{\star sp}$. Then, $t\in O^{\star s}$ and $t\in O^{\star p}$. This implies, $O\cap U\notin \mathbb{I}_{T}$ and $O\cap V\notin \mathbb{I}_{T}$ for all $U\in \mathcal{O}^s (t)$  and $V\in \mathcal{O}^p (t)$. Thus, $O\cap U\neq \emptyset$ and $O\cap V\neq \emptyset$ for all $U\in \mathcal{O}^s (t)$  and $V\in \mathcal{O}^p (t)$ and hence $t\in sCl(O)$ and $t\in pCl(O)$. Thus, $t\in sCl(O)\cap pCl(O)$ and Consequently $O^{\star sp}\subseteq sCl(O)\cap pCl(O)$.

Conversely, let $q\in sCl(O)\cap pCl(O)$. Then, $q\in sCl(O)$ and $q\in pCl(O)$ and hence $O\cap U\neq \emptyset$ and $O\cap V\neq \emptyset$ for all $U\in \mathcal{O}^s (q)$  and $V\in \mathcal{O}^p (q)$. Thus, $O\cap U\notin \mathbb{I}_{T}$ and $O\cap V\notin \mathbb{I}_{T}$ for all $U\in \mathcal{O}^s (q)$  and $V\in \mathcal{O}^p (q)$ as $\mathbb{I}_{T} \cap \mathcal{O}^s = \{\emptyset \}$ and $\mathbb{I}_{T} \cap \mathcal{O}^p = \{\emptyset \}$. This implies, $q\in O^{\star s}$ and $q\in O^{\star p}$ and hence  $q\in O^{\star sp}$. Thus, $sCl(O)\cap pCl(O)\subseteq O^{\star sp}$.

This completes the proof.
\end{proof}
\begin{prop}
\label{R5.1.12}
Let $\mathbb{I}_{T}$ be an ideal on a topological space $\mathbb{T}$. Then, $\mathbb{I}_{T} \cap \mathcal{O}^s = \{\emptyset\}$ and $\mathbb{I}_{T} \cap \mathcal{O}^p = \{\emptyset\}$ if and only if $\mathbb{T}=\mathbb{T}^{\star sp}$.
\end{prop}
\begin{proof}
Let $T={T}^{\star sp}$ and $\emptyset \neq U\in \mathcal{O}^s$  and $\emptyset \neq V\in \mathcal{O}^p$. Then, $U=U\cap {T}\notin  \mathbb{I}_{T}$ and $V=V\cap {T}\notin  \mathbb{I}_{T}$. Since, $U\in \mathcal{O}^s$  and $V\in \mathcal{O}^p$ are arbitrary, then $\mathbb{I}_{T} \cap \mathcal{O}^s = \{\emptyset \}$ and $\mathbb{I}_{T} \cap \mathcal{O}^p = \{\emptyset \}$.

Conversely, assume that $\mathbb{I}_{T} \cap \mathcal{O}^s = \{\emptyset \}$ and $\mathbb{I}_{T} \cap \mathcal{O}^p = \{\emptyset \}$ hold. We have to show that ${T}\setminus {T}^{\star sp}$ must be empty. If not, assume there exists $t\in {T}$ such that $t\notin {T}^{\star sp}$. This implies, there exists either $\emptyset \neq U\in \mathcal{O}^s (t)$  or $\emptyset \neq V\in \mathcal{O}^p (t)$ such that either $U=U\cap {T}\in  \mathbb{I}_{T}$ or $V=V\cap {T}\in  \mathbb{I}_{T}$ which contradicts our given condition.

This completes the proof.
\end{proof}

\begin{Defn}
\label{ }
Let $\mathbb{I}_T$ be an ideal on a topological space $\mathbb{T}$ and $S\subseteq {T}$. Then, $S$  is called $\star sp$-dense in $\mathbb{T}$ if $S^{\star sp}={T}$.

Equivalently, $S$ is $\star sp$-dense in $\mathbb{T}$ iff $S\cap U\notin \mathbb{I}_{T}$ and $S\cap V\notin \mathbb{I}_{T}$ for all $\emptyset\neq U\in \mathcal{O}^s $  and $\emptyset\neq V\in \mathcal{O}^p$.
\end{Defn}

For the existence of $\star sp$-dense set, we consider the following example:
\begin{Exmp}
\label{Example 5.1.2}
In Example \ref{Example 5.1.1}, if we consider $S=\{t_1, t_2\}$, then $S^{\star sp}={T}$. Hence, $S$ is a $\star sp$-dense subset of ${T}$.
\end{Exmp}

\begin{prop}
\label{R5.1.2}
Let $\mathbb{I}_{T}$ be an ideal on a topological space $\mathbb{T}$. Then, every $\star sp$-dense subset of ${T}$ is $\mathbb{I}_{T}$-dense subset in $\mathbb{T}$.
\end{prop}
\begin{proof}
Let $S$ be a $\star sp$-dense subset of ${T}$. Then, for any  $t\in {T}$,   $S\cap U\notin \mathbb{I}_{T}$ and $S\cap V\notin \mathbb{I}_{T}$ for all $U\in \mathcal{O}^s (t)$  and $V\in \mathcal{O}^p (t)$. Since, $\mathcal{O}\subseteq\mathcal{O}^s$ and $\mathcal{O}\subseteq\mathcal{O}^p$, then $S\cap O\notin \mathbb{I}_{T}$ holds for any $O\in \mathcal{O}(t)$. Hence,  $t\in S^{\star}$. Since $t$ is an arbitrary, ${T}\subseteq S^{\star}$. This completes the proof.
\end{proof}
The reverse inclusion may not be true. For this, we consider the following example:
\begin{Exmp}
Consider an ideal topological space  $({T}, \tau_{T}, \mathbb{I}_T)$ where $T =\{t_1, t_2, t_3 \},\; \tau_{T} = \{\emptyset, T, \{t_1, t_2 \} \},\; $ $\mathbb{I}_{T} = \{\emptyset, \{t_1\}, \{t_3 \}, \{t_1, t_3 \}\}$. Then $SO(\mathbb{T}) = \{\emptyset, T, \{t_1, t_2 \} \},\;  PO(\mathbb{T}) = \{\emptyset, T, \{t_1, t_2 \}, \{t_1 \}, \{ t_2 \}, \{t_1, t_3 \}, $ $ \{t_2, t_3 \} \}$. Take $A = \{t_1, t_2 \}$, then $A^{\star} = T$, $A^{\star p} = \{t_2\}$ and  $A^{\star s} = T$. So,  $A^{\star sp}=\{t_2\} $. Hence, $A$ is $\mathbb{I}_{T}$-dense subset in $\mathbb{T}$ but not $\star sp$-dense.
\end{Exmp}
Clearly if $\mathbb{I}_{T} \cap \mathcal{O}^s \neq \{\emptyset\}$ and $\mathbb{I}_{T} \cap \mathcal{O}^{p} \neq \{\emptyset\}$, then no subset of $T$ is ${\star sp}$-dense not even ${T}$ itself.
\begin{Defn}
\label{5.2}
An ideal  topological space $(T, \tau_{T}, \mathbb{I}_{T})$  is called $\star sp$-resolvable if and only if $T$ has two disjoint $\star sp$-dense subsets.
\end{Defn}

\begin{prop}
\label{R5.1.3}
If an ideal  topological space $(T, \tau_{T}, \mathbb{I}_{T})$  is $\star sp$-resolvable, then $\mathbb{I}_{T} \cap \mathcal{O}^s = \{\emptyset\}$ and $\mathbb{I}_{T} \cap \mathcal{O}^p = \{\emptyset\}$.
\end{prop}
\begin{proof}
Let $\bf A$ be subset $(T, \tau_{T}, \mathbb{I}_{T})$ such that $\bf A$ forms $\star sp$-resolution with its complement. It follows that $T$ is $\star sp$-dense. Then, by Result \ref{R5.1.12}, $\mathbb{I} \cap \mathcal{O}^s = \{\emptyset\}$ and $\mathbb{I} \cap \mathcal{O}^p = \{\emptyset\}$ hold.
\end{proof}

\begin{Defn}
\label{5.1}
Let $\mathbb{I}_{T}$ be an ideal on a topological space $\mathbb{T}$ and $\zeta : T\longmapsto T$ be a mapping. Then the dynamical system $(\mathbb{T}, \zeta)$ is called ${\star sp}$ transitive if for each pair of nonempty open sets $U, V$, there exists $n\in \mathbb{N}$ such that $(\zeta^n(U))^{\star sp} \cap V^{\star sp}\neq \emptyset$.
\end{Defn}
We are now considering the following example for the existence of Definition \ref{5.1}:
\begin{Exmp}
\label{Example 5.1}
Suppose $T=\{t_1, t_2, t_3\}$ is a set endowed with the topology $\tau_{T}=\{\emptyset, T, \{t_1\}, \{t_2\}, \{t_1, t_2\}\}$ and $\mathbb{I}_{T}=\{\emptyset,\{t_3\}\}$. Now $PO(\mathbb{T})= \{\emptyset, T, \{t_1\}, \{t_2\}, \{t_1, t_2\}\}$ and $SO(\mathbb{T})= \{\emptyset, T, \{t_1\}, \{t_2\}, \{t_1, t_2\}, $ $ \{t_2, t_3\}\}$. Now, $\{t_1\}^{*s}=\{t_1\}$, $\{t_1, t_2\}^{*s}=T$, $\{t_2, t_3\}^{*s}=\{t_2, t_3\}$, $\{t_1, t_3\}^{*s}=\{t_1, t_3\}$, $\{t_2\}^{*s}=\{t_2, t_3\}$, $\{t_3\}^{*s}=\{t_3\}$, $\{\emptyset\}^{*s}=\{\emptyset\}$, $\{T\}^{*s}={T}$ and $\{t_1\}^{*p}=\{t_1, t_3\}$, $\{t_2\}^{*p}=\{t_2, t_3\}$, $\{t_3\}^{*p}=\emptyset$, $\{t_1, t_2\}^{*p}={T}$, $\{t_2, t_3\}^{*p}=\{t_2, t_3\}$, $\{t_3, t_1\}^{*p}=\{t_3, t_1\}$, $\{\emptyset\}^{*p}=\{\emptyset\}$, $\{T\}^{*p}={T}$.  Let us define a mapping $\zeta:{T} \rightarrow {T}$ by $\zeta(t_1)=t_2$, $\zeta(t_2)=t_1$ and  $\zeta(t_3)=t_3$. Then,  $(\zeta^n(U))^{\star sp} \cap V^{\star sp}\neq \emptyset$, for each pair of nonempty open sets $U, V$ for all odd positive integers $n$. Thus, the dynamical system $(T, \zeta)$ is ${\star sp}$-transitive.
\end{Exmp}

\begin{note}
\label{N5.1}
Let $\mathbb{I}_{T}$ be an ideal on a topological space $\mathbb{T}$ and $\zeta : \mathbb{T}\longmapsto \mathbb{T}$ be a mapping. If the dynamical system $(\mathbb{T}, \zeta)$ is  ${\star sp}$ transitive, then for each  nonempty open set $U$ in $\mathbb{T}$, $\bigcup\limits_{n=0}^{\infty}f^{n}(U)$  may not be ${\star sp}$-dense in $\mathbb{T}$.
\end{note}
We are now giving an example in support of the Note \ref{N5.1}:
\begin{Exmp}
Consider an ideal topological space  $(T, \tau_{\mathbb{T}}, \mathbb{I}_T)$ where $T =\{t_1, t_2, t_3 \},\; \tau_{T} = \{\emptyset, T, \{t_1, t_2 \} \},\;  \mathbb{I}_{T} = \{\emptyset, \{t_1\}, \{t_3 \}, \{t_1, t_3 \}\}$. Then $SO(\mathbb{T}) = \{\emptyset, T, \{t_1, t_2 \} \},\;  PO(\mathbb{T}) = \{\emptyset, T, \{t_1, t_2 \}, \{t_1 \}, \{ t_2 \}, \{t_1, t_3 \}, $ $ \{t_2, t_3 \} \}$. Take $U = \{t_1, t_2 \}$, then $U^{\star p} = \{t_2\}$ and  $U^{\star s} = T$. So,  $U^{\star sp}=\{t_2\} $. Let us define a mapping $\zeta:T \rightarrow T$ by $\zeta(t_1)=t_2$, $\zeta(t_2)=t_1$ and  $\zeta(t_3)=t_3$. Then,  $(\zeta^n(U))^{\star sp} \cap V^{\star sp}\neq \emptyset$, for each pair of nonempty open sets $U, V$ for all  positive integers $n$. Thus, the dynamical system $(T, \zeta)$ is ${\star sp}$-transitive. Now, $(\bigcup\limits_{n=0}^{\infty}f^{n}(U))^{\star sp}=\{t_2\}$ and hence $\bigcup\limits_{n=0}^{\infty}f^{n}(U)$  is not ${\star sp}$-dense in $T$
\end{Exmp}
\begin{Defn}
\label{5.2}
Let $\mathbb{I}_T$ be an ideal on a topological space $\mathbb{T}$ and $\zeta : T \longmapsto T$ be a mapping. A point $t\in T$ is called ${\star sp}$ non-wandering point if for each nonempty open set $U$ containing $t$, there exists $n\in \mathbb{N}$ such that $(\zeta^n(U))^{\star sp} \cap U^{\star sp}\neq \emptyset$.
\end{Defn}
\begin{Exmp}
\label{Example 5.2}
Suppose $T=\{t_1, t_2, t_3\}$ is a set endowed with the topology $\tau_{T}=\{\emptyset, T, \{t_1\}, \{t_2\}, \{t_1, t_2\}\}$ and $\mathbb{I}_{T}=\{\emptyset,\{t_3\}\}$. Now $PO(\mathbb{T}) = \{\emptyset, T, \{t_1\}, \{t_2\}, \{t_1, t_2\}\}$ and $SO(\mathbb{T})= \{\emptyset, T, \{t_1\}, \{t_2\}, \{t_1, t_2\}, $ $ \{t_2, t_3\}\}$. Now, $\{t_1\}^{*s}=\{t_1\}$, $\{t_1, t_2\}^{*s}=T$, $\{t_2, t_3\}^{*s}=\{t_2, t_3\}$, $\{t_1, t_3\}^{*s}=\{t_1, t_3\}$, $\{t_2\}^{*s}=\{t_2, t_3\}$, $\{t_3\}^{*s}=\{t_3\}$, $\{\emptyset\}^{*s}=\{\emptyset\}$, $\{T\}^{*s}=T$ and $\{t_1\}^{*p}=\{t_1, t_3\}$, $\{t_2\}^{*p}=\{t_2, t_3\}$, $\{t_3\}^{*p}=\emptyset$, $\{t_1, t_2\}^{*p}=T$, $\{t_2, t_3\}^{*p}=\{t_2, t_3\}$, $\{t_3, t_1\}^{*p}=\{t_3, t_1\}$, $\{\emptyset\}^{*p}=\{\emptyset\}$, $\{T\}^{*p}=T$.  Let us define a mapping $\zeta:T \rightarrow T$ by $\zeta(t_1)=t_2$, $\zeta(t_2)=t_1$ and  $\zeta(t_3)=t_3$. Take, $t_1\in \mathbb{T}$. The open sets containing $t_1$ are $T, \{t_1\}, \{t_1, t_2\}$. Then,  $(\zeta^n(U))^{\star sp} \cap U^{\star sp}\neq \emptyset$, for each nonempty open set $U$ containing $t_1$ when $n=2$. Thus, $t_1$ is ${\star sp}$ non-wandering point of $T$.
\end{Exmp}

\begin{thm}
\label{ }
Let $\mathbb{I}_T$ be an ideal on a topological space $\mathbb{T}$ and $(\mathbb{T}, \zeta)$ is a dynamical system. If $(\mathbb{T}, \zeta)$ is ${\star sp}$ transitive then for all $(\emptyset\neq) O \in \mathcal{O}$, $O \notin \mathbb{I}_T$ as well as $\zeta^{n}(O) \notin \mathbb{I}_T$ for some positive integer $n$.
\end{thm}
\begin{proof}
  The proof is obvious and hence omitted.
\end{proof}
\begin{thm}
\label{ }
Let $\mathbb{I}_T$ be an ideal on a topological space $\mathbb{T}$. Suppose $O(\neq \emptyset) \in \mathcal{O},\; U(\neq \emptyset) \in \mathcal{O}^s,\; V(\neq \emptyset) \in \mathcal{O}^p$. Then for any $t \in O,\; U$ or $V$, $t \notin (\zeta^n(W))^{\star sp}$ for every $W \in \mathcal{O}$.
\end{thm}

\section{Conclusion}

The local functions in a topological space have been extended to the generalized local functions through this article. As a result of this research, multiple branches have emerged.  In this work, we use generalized open sets, such as semi-open sets, preopen sets, $b$-open sets, and $\beta$-open sets, to study the local function in greater detail and analyse its various facets. In this article, a new kind of dense set in topological spaces has been developed, and its relationship among generalized open sets has been discussed.  For the apllication of the $\star sp$ local function, we have also investigated the dynamical system and topologically transitive scenario.


\providecommand{\bysame}{\leavevmode\hbox
to3em{\hrulefill}\thinspace}


\end{document}